\documentclass[11pt,notitlepage,twoside,a4paper]{amsart}
 \usepackage{amsfonts}
 \usepackage{mathrsfs}
\usepackage{amsmath,amssymb,enumerate}

\usepackage{epsfig,fancyhdr,color}

\usepackage{amssymb}
\usepackage{amsmath,amsthm}
\usepackage{latexsym}
\usepackage{amscd}
\usepackage{psfrag}
\usepackage{graphicx}
\usepackage[latin1]{inputenc}
\usepackage[all]{xy}

\newtheorem{theorem}{\rm\bf Theorem}[section]

\newtheorem*{theorem 1}{\rm\bf Proposition 1}
\newtheorem*{theorem 2}{\rm\bf Proposition 2}
\newtheorem*{proposition 1}{\rm\bf Proposition 3.2}

\theoremstyle{definition}
\newtheorem{definition}[theorem]{\rm\bf Definition}

\theoremstyle{remark}
\newtheorem{remark}[theorem]{\rm\bf Remark}

\newtheorem{example}[theorem]{\rm\bf Example}

\def\interieur#1{\mathord{\mathop{\kern 0pt #1}\limits^\circ}}

\definecolor{NoteColor}{rgb}{1,0,0}

\title[Arc length for locally analytic curves]{Arc length as a conformal parameter for locally analytic curves}
\author{Vassili Nestoridis and Athanase Papadopoulos}
\address{V. Nestoridis, Department of Mathematics, University of Athens, 15784 Panepistemioupolis, Athens, Greece,
 email: 
 {\rm vnestor@math.uoa.gr}}
 \address{A. Papadopoulos,
Institut de Recherche Math\'ematique Avanc\'ee,
(Universit\'e de Strasbourg et CNRS),
7 rue Ren\'e Descartes,
67084 Strasbourg Cedex France,
 email: {\rm papadop@math.unistra.fr} }

\date{\today}

\begin{document}

\maketitle

\centerline{{\it Dedicated to Professor Richard Aron on the occasion of his retirement
}}


\begin{abstract}
For any locally analytic curve we show that arc length can be complexified and seen as a conformal parameter. As an application, we show that any such curve defines a unique maximal one and that the notions of analytic Jordan curve coincides with the notion of a Jordan curve which is locally analytic. We give examples where we also find, for a curve $\gamma(s)$, the limit sets of the largest extensions, that is,  the limit set of the curve as $s$ converges to the end point of the interval. \end{abstract}

\noindent AMS Mathematics Subject Classification: 30B40, 32B15
\medskip

\noindent Keywords: Analytic curve, locally analytic curve, conformal parameter, arc length, analytic extension.
\medskip
\maketitle

%
%
%
%
%

\section{Introduction}

Let us first recall the notion of analytic (compact) curve and of analytic Jordan curve, following \cite{Ahlfors}. These are respectively the images of $X=[0,1]$ or $X=\mathbb{T}=\{z\in\mathbb{C}:\vert z\vert =1\}$ (the unit circle) by conformal mappings defined on an open set $V$ such that $X\subset V\subset \mathbb{C}$. We naturally extend the notion of analytic curve to the notion of locally analytic curve. Such a curve may be either simple or not (that is, with auto-intersection). A natural question is whether the notion of analytic Jordan curve coincides with the notion of locally analytic Jordan curve. We show that the answer is positive. In order to prove this equivalence, we are led to parametrize an arbitrary analytic curve by arc length and prove that arc length can be extended to take complex values and become a conformal parameter. This implies that all the extensions of a locally analytic curve are compatible and define a maximal analytic curve which is unique. More generally, if $\gamma$ is a map defined on an open interval $I$ of $\mathbb{R}$
with complex values, with non-vanishing derivative and locally
analytic, then following this parametrization of the curve we obtain a maximal one having the same properties and extending $\gamma$ to the largest possible interval $\tilde{I}$ in $\mathbb{R}$ containing $I$. Our proof shows that using arc length, we get at least the same extension.  Example \ref{ex:5} below shows that in some cases the extension with arc length is larger. In this sense, the best conformal parametrization of any locally analytic curve is by arc length. 

We then give examples of such curves.
The limit sets in these examples are singletons $\{w_0\}$, $w_0\in\mathbb{C}\cup\{\infty\}$, or a circle. 

We address the question of characterizing such limit sets. A known characterization of analytic Jordan curves is that they are of the form
\[\phi(\{z\in\mathbb{C}:\vert z\vert =r\})\ 0<r<1,\]
where $\phi$ is a conformal map defined on the unit disc. Our characterization of locally analytic curves is that they are of the form $\phi(I)$, where $I\subset \mathbb{R}$ is an open interval and $\phi$ is holomorphic and locally injective on an open set $V$, $I\subset V\subset \mathbb{C}$. Note that if we consider the restriction of $\phi$ on a (not necessarily open) subinterval $I'$ of $I$  we also obtain a locally analytic curve. By the definition of a locally analytic curve the converse also holds, that is, we have a characterization.

Finally, a corollary of the previous result is that if $\gamma=\phi(e^{i\theta})$ is an analytic Jordan curve and $f$ a complex function defined on the image of $\gamma$, then saying that $f$ is analytic at any point of $\gamma$ with respect to the variable $\theta$ is equivalent to saying that it is analytic with respect to the arc length $s$. 

\section{The results}

We start with a classical definition from Ahlfors' book \cite{Ahlfors}:
\begin{definition}\label{curve}
Let $a<b$ be two real numbers and $\gamma:[a,b]\to\mathbb{C}$ a continuous mapping. The curve $\gamma$ is said to be analytic if there exists an open set $V$ in $\mathbb{C}$ containing the segment $[a,b]$ and a conformal mapping $\phi: V\to \mathbb{C}$ such that $\phi\vert_{[a,b]}=\gamma$. 
\end{definition}
\begin{remark}
It is easy to see that in Definition \ref{curve} one may assume without loss of generailty that $V$ is an open simply connected domain containing $[a,b]$ and symmetric with respect to the real axis, that is, $\overline{z}\in V$ for any $z$ in $V$.
Combining with the result of \cite{2},   we see that it suffices that $\phi$ is holomorphic in some open set $W\subset \mathbb{C}$ containing $[a,b]$, that $\phi\vert_{[a,b]}$ is 1-1 and that $\phi'(z)\not=0$ for each $z\in[a,b]$. With these conditions, we can find a simply connected domain $V$ containing $[a,b]$ symmetric with respect to the real axis such that $\phi\vert_V$ is conformal (and $\phi\vert_{[a,b]}=\gamma$); hence the formulation in Definition \ref{curve}.
\end{remark}

\begin{definition}\label{Jordan}
A Jordan curve $\gamma$ in $\mathbb{C}$ is said to be an analytic Jordan curve if it is the image $F(\mathbb{T})$ of the unit circle by a conformal mapping $F$ defined on an open annulus $\mathrm{D}$ containing $\mathbb{T}$. In other words, there exists $F:\mathrm{D}\to\mathbb{C}$ which is 1-1 holomorphic, with 
\[\mathrm{D}=\{z\in\mathbb{C}:r_1< z<r_2\},\ 0<r_1<1<r_2<+\infty\]
and $\gamma=F(\mathbb{T})$.
\end{definition}

We extend Definition \ref{curve} by considering curves which locally satisfy that definition.

\begin{definition}\label{curve-loc}
Let $I\subset\mathbb{R}$ be an interval of any type and $\gamma:I\to\mathbb{C}$ a continuous function. Let $t_0\in I$. We say that $\gamma$ satisfies locally at $t_0$ Definition \ref{curve} if the following holds:
\begin{enumerate}
\item For any $t_0\in \overset{o}I$, there exist $a,b\in I$ with $a<t_0<b$ such that $\gamma\vert_{[a,b]}$ satisfies Definition \ref{curve}.
\item If $t_0=\sup I=\max I \in I\subset \mathbb{R}$, then there exists $a\in I$ with $a<t_0$ such that $\gamma\vert_{[a,t_{0}]}$ satisfies Definition \ref{curve} with $b=t_0$.
\item If $t_0\in\mathbb{R}$ and $t_0=\inf I=\min I \in I\subset \mathbb{R}$, then there exists $b\in I$ with $t_0< b$ such that $\gamma\vert_{[t_{0},b]}$ satisfies Definition \ref{curve} with $a=t_0$.
\end{enumerate}
Finally, we say that $\gamma$ satisfies locally Definition \ref{curve} (or that $\gamma$ is locally an analytic curve) if it satisfies Definition \ref{curve} at any point $t_0$ of its domain of definition $I$.
\end{definition}
We note that if $I=(a,b]$ or $I=[a,b)$ or $I=(a,b)$ and $\gamma$ satisfies Definition \ref{curve-loc} on $I$ then there exists an open interval $I'$, $I\subset I'\subset \mathbb{R}$ such that $\gamma$ extends to $I'$ and the extension satisfies Definition \ref{curve-loc} on $I'$.

Therefore, we can restrict our attention to the case where $I$ is an open subinterval of $\mathbb{R}$.

It is easy to see that every analytic Jordan curve (Definition \ref{Jordan}) satisfies locally Definition \ref{curve}, that is, it  satisfies Definition \ref{curve-loc} on its domain of definition $I=\mathbb{R}$, where it is a periodic function. A natural question is whether the converse holds. That is, if a Jordan curve locally satisfies Definition \ref{curve}, is it true that it is an analytic Jordan curve and satisfies Definition \ref{Jordan}?

In this article, we give an affirmative answer to this question, by showing that if two curves locally satisfy Definition \ref{curve} and if they have a common arc of strictly positive length, then one extends the other and their union defines a curve locally satisfying Definition \ref{curve}. In other words, if a curve locally satisfies Definition \ref{curve} and can be extended and still has the same property, then it has only one possible (maximal) extension.

The idea of the proof is to use arc length parametrization. Suppose that $\gamma$ satisfies Definition \ref{curve}. The length $s=s(t)$, $t\in [a,b]$ is a real parameter satisfying 
\[\frac{d s}{dt}=\vert \gamma'(t)\vert=  \vert \phi'(t)\vert= \sqrt{\phi'(t)\overline{\phi'(t)}}=\sqrt{\phi'(t)\overline{\phi'(\overline{t})}}\]
on $[a,b]$. The question is then whether one can extend this to the complex domain. That is, whether the function $s:[a,b]\to\mathbb{R}$ can be extended to a holomorphic function on an open set $V$, $[a,b]\subset V\subset \mathbb{C}$. Then $s(z)$ will be complex, and the variable $z$ as well. If in addition the extension of $s$ is conformal on $V$, then it can be seen as a conformal parameter and it can replace $\phi$ in Definition \ref{curve}. Then, if $L>0$ is the length of $\gamma$ on $[a,b]$ and $s_0\in\mathbb{R}$, there exists an open set $\mathcal{O}$, $[s_0,s_{0+L}]\subset \mathcal{O}\subset\mathbb{C}$ and a conformal mapping $\delta:\mathcal{O}\to \mathbb{C}$ such that $\gamma=\delta\circ s$, $\gamma=\delta\vert_{[s_{0},s_{0}+L]}$, that is, $\gamma(t)=\delta(s(t))$ for $t\in[a,b]$, and this representation is unique modulo real translations. Of course, $\mathcal{O}$ can be chosen to be a simply connected domain symmetric with respect to the real axis;  in fact  $\mathcal{O}$  can be chosen to be a rectangle
symmetric with respect to the real axis.

If there are two extensions of $\delta$ around $s_0+L$, they will be holomorphic and they will coincide on $(s_0,s_0+L)$. By analytic continuation, they will coincide on an open set containing $(s_0,s_0+L+\eta)$ for some $\eta>0$. Thus, the extension, if it exists, is unique.

In order to justify the existence of the extension of the map $s:[a,b]\to \mathbb{C}$ to an open set $V$, $[a,b]\subset V\subset \mathbb{C}$, we consider the function $\phi'(z)\overline{\phi'(\overline{z})}$ which is holomorphic on $V$ and does not vanish. Since $V$ is simply connected, there exists a holomorphic branch of $\sqrt{\phi'(z)\overline{\phi'(\overline{z})}}$ on $V$. We choose it so that on $[a,b]$ it takes positive values and for $z=a$ it takes the value $\vert \phi'(a)\vert >0$.

Since $V$ is simply connected, every holomorphic function on $V$ has a primitive on $V$. Such a primitive $G$ of $\sqrt{\phi'(z)\overline{\phi'(\overline{z})}}$ on $V$ is the complex holomorphic extension of the map $s:[a,b]\to\mathbb{C}$. Since $G\vert_{[a,b]}=s$ is 1-1 and $G'(t)=s'(t)\not=0$ on $[a,b]$ because $s'(t)=\vert\gamma'(t)\vert=\vert\phi'(t)\vert>0$, using the result of \cite{2}, we see that $G$ is conformal on an open set containing $[a,b]$, and this open set can certainly be chosen to be simply connected and symmetric with respect to the real axis. Then, it suffices to set $\delta=G$. Continuing in this way, we obtain a locally analytic maximal extension.

Thus, we have the following:
\begin{theorem}
Let $\gamma$ be a curve satisfying locally Definition \ref{curve}. Then, the union of the extensions of $\gamma$ which are locally analytic define a maximal curve $\gamma^*:(A,B)\to \mathbb{C}$, $-\infty\leq A<B\leq+\infty$, parametrized by arc length, which is locally analytic. The mapping $\gamma^*$ has a holomorphic extension to an open set in $\mathbb{C}$ containing $(A,B)$ which is locally conformal (and locally injective).
\end{theorem}

We note that although $\gamma^*$ is locally injective, it is not necessarily injective on its domain of definition $(A,B)$. A case where $\gamma^*$ is globally injective on $(A,B)=\mathbb{R}$ is when $\gamma$ is a straight line. If $\gamma$ is an analytic Jordan curve, then $A=-\infty, B=+\infty$ and $\gamma^*$ is periodic, with minimum positive period the length $L>0$ of $\gamma$. By analytic continuation the locally conformal extension $G$ of $\gamma^*$ on an open domain containing $\mathbb{R}$ also satisfies $G(z+L)=G(z)$ for all $z$. (It is easy to see that the domain of definition of $G$ contains a band $\{z\in\mathbb{C}: \vert \mathrm{Im}(z)\vert <\eta\}$ for some $\eta>0$.) Without loss of generality we assume $L=2\pi$. Then, for $w=re^{i\theta}$ on an open annulus $\mathcal{D}$ containing the unit circle $\mathbb{T}$, we set 
\[\phi(w)=G\left(\frac{\ln w}{i}\right).\] This function is well defined on $\mathcal{D}$. Since $\phi(e^{i\theta})\not=0$ and $\phi\vert_{\mathbb{T}}$ is 1-1, then by \cite{2} we conclude that $\phi$ is 1-1 and conformal on an open annulus $\mathcal{D}'$,  $\mathbb{T}\subset \mathcal{D}'\subset \mathcal{D}$. 

Thus, we have the following
\begin{theorem}
Let $\gamma$ be Jordan curve which is a locally analytic curve in the sense of Definition \ref{curve-loc}. Then $\gamma$ is an analytic Jordan curve satisfying Definition \ref{Jordan}.
\end{theorem}
\begin{remark}
If $\gamma$ is a locally analytic closed curve (not necessarily Jordan), then $A=-\infty$, $B=+\infty$ and $\gamma^*$ has a holomorphic extension to an open set containing $\mathbb{R}$ which is locally injective, conformal and periodic. The  converse also holds. For an example of such a curve which is not Jordan, see Example \ref{ex:7} below. 
\end{remark}
\begin{remark}
The map $\gamma^*:(A,B)\to\mathbb{C}$ is a canonical representation of the maximal extension of a locally analytic curve, which is unique modulo real translations. The fact that $\gamma^*$ can be extended to an open set in $\mathbb{C}$ containing $(A,B)$ and the extension is holomorphic locally injective, is equivalent to saying that the derivative of $\gamma^*$ does not vanish at any point of $(A,B)$ and that $\gamma^*$ is locally analytic; that is, for every $t_0\in(A,B)$, there exists $\eta>0$, $\eta \leq \min(B-t_0,t_0-A)$ and a power series $\sum_{n=0}^\infty a_n (t-t_0)^n$, $a_n\in\mathbb {C}$ converging on 
$(t_0-\eta,t_0+\eta)$ such that $\gamma^*(t)=\sum_{n=0}^\infty a_n(t-t_0)^n$ on $(t_0-\eta,t_0+\eta)$ and $a_1\not=0$. Therefore, the definition of locally analytic curve (Definition \ref{curve-loc}) could be changed to be that $\gamma:(a,b)\to\mathbb{C}$ is locally analytic if (and only if) $\gamma'(t)\not=0$ for all $t\in(a,b)$ and $\gamma$ is real analytic on $(a,b)$. If the domain of definition is not an open interval, then we require that there is an extension of the type described. 

The result that we proved says that the above parameter $t$ can always be chosen to coincide with the arc length of $\gamma$, that we have real analyticity with respect to $s$, and that the derivative with respect to $s$ does not vanish at any point.
\end{remark}

\begin{remark}
Let $\gamma$ be an analytic Jordan curve and $\phi:\mathbb{D}\to\mathbb{C}$ a conformal mapping from the open unit disc $\mathbb{D}$ onto the interior of $\gamma$. Then $\phi$ has a conformal extension on a disc $D_r=\{z\in\mathbb{C}:\vert z\vert<r\}$ with $r>1$ (see \cite{Ahlfors}).
The curve $\gamma$ is the curve $\phi(e^{it})$, $t\in\mathbb{R}$. Let $s:\mathbb{R}\to\mathbb{R}$ be the arc length function of $\gamma$. Then by the previous results $s$ is invertible: $s^{-1}:\mathbb{R}\to\mathbb{R}$, $\tau=s(t)$, $t=s^{-1}(\tau)$ and both functions have non-vanishing derivatives and they are analytic at every point. It follows that for a function $f:\gamma\to \mathbb{C}$ saying that $f$ is analytic at any point of $\gamma$ with respect to the parameter $t$ is equivalent to being analytic with respect to the arc length of $\gamma$ at any point of $\gamma$.

Equivalently, $f$ can be holomorphically extended on an open set $V:\gamma\subset V\subset\mathbb{C}$.

\end{remark}

\begin{remark}
Our result implies that if a curve is locally analytic with
respect to different (but compatible)  parametrizations, then there
exists one parametrization for the whole curve with respect to which the
curve becomes locally analytic. In case the curve is injective (that is, a
simple curve), then the arc length is a function of the position and
according to our result it is a conformal parameter for the whole curve.
In the general case, locally, we have an injective curve and this allows us
to continue and cover the whole curve.
\end{remark}

\section{Examples}
As we already mentioned, the map $\gamma^*:(A,B)\to \mathbb{C}$ is not injective in general; it is only locally injective. If the curve $\gamma$ has infinite length, then $A=-\infty$ or $B=+\infty$. In Example \ref{ex:2} below we see that the length of $\gamma$ is finite but since $\gamma^*$ extends by periodicity we have $A=-\infty$ and $B=+\infty$. Now we give some examples.

\begin{example}\label{ex:1}

Let $\gamma$ be the real line; then $A=-\infty, B=+\infty$ and $\gamma^*(s)=s+s_0$, $s\in \mathbb{C}$.
\end{example}

\begin{example}\label{ex:2}
Let $\gamma$ be the unit circle; then $A=-\infty, B=+\infty$ and $\gamma^*(s)=e^{is}$, $s\in \mathbb{C}$. 
\end{example}

\begin{example}\label{ex:3}
Let $\gamma(t)=\displaystyle (e^{it}-1)\exp\frac{e^{it}+1}{e^{it}-1}$, $0<t<2\pi$. Then $\gamma$ is 	a double spiral pointing from $0$ to the same point $0$, it has infinite length, $A=-\infty$ and $B=+\infty$. It satisfies Definition \ref{curve-loc} and it is locally injective but not globally. We also have
\[\lim_{t\to 0^{+}}\gamma(t)=\lim_{\tau\to -\infty}\gamma^*(\tau)
=\lim_{t\to 2\pi^{-}}\gamma(t)=\lim_{\tau\to\infty}\gamma^*(\tau)=0.\]
\end{example}

\begin{example}\label{ex:4}
Let $\gamma(x)=e^xe^{+i/x}$, $-\infty<x<0$. This map extends conformally on an open set containing $(-\infty,0)$. The length of $\gamma$ is infinite as $x\to 0^-$; thus $B=+\infty$. The limit set of $\gamma$ as $x\to 0^-$ (equivalently, $s\to+\infty$) is the unit circle $\mathbb{T}$, since $\gamma(x)=e^xe^{+i/x}$ and $e^x\to 1$ while $e^{+i/x}$ covers $\mathbb{T}$ infinitely many times as $x\to0^-$.
As $x\to-\infty$ we have $\gamma(x)\to 0$, the positive real axis is a tangent of $\gamma$ at $0=\gamma(-\infty)$ and the length is finite.
We should investigate whether the curve $\gamma$ can be continued beyond this point while it still  satisfies Definition \ref{curve-loc}.

Let us set $y=-x\to+\infty$. 
Then 
\[\delta(y)=\gamma(x)=e^{-y-i/y}=\frac{1}{e^{y+i/y}}\]
and we are interested in the case $y\to+\infty$.

Then, $\lim_{y\to+\infty} \delta(y)=0$. Therefore, it is possible that with another parametrization, the curve can be further extended while it stays locally analytic. Then, this will hold for the parameter arc length. Thus, let us compute the derivative $d\delta /dy$, whose absolute value is the derivative of arc length.

\[\frac{d\delta}{dy}=\delta'(y)=e^{-y-i/y}\cdot (-1+\frac{i}{y^2})=\frac{-y^2+i}{e^{y+i/y}y^2},
\]
$\vert\delta'(y)\vert = \frac{\sqrt{1+y^4}}{e^yy^2}$ which is integrable as $y\to+\infty$.

If we can extend the curve, then $ds/d\delta$ and its derivatives will have finite limit as $y\to+\infty$. Let us compute these derivatives.

\[\frac{ds}{d\delta}=\frac{ds}{dy}\cdot \frac{1}{d\delta /dy}=\vert \delta'(y)\vert\cdot \frac{1}{\delta'(y)}
\]
\[= \frac{\sqrt{1+y^4}}{e^yy^2}\frac{y^2e^{y+i/y}}{-y^2+i}
=\frac{\sqrt{1+y^4}}{-y^2+i}e^{i/y}
.\]
We let $\phi(y)$ denote the last term. We have
\[\frac{d^2s}{d\delta^2}=\frac{d\phi}{d\delta}=\frac{d\phi}{dy}\cdot \frac{1}{d\delta /dy}=\phi'(y)\frac{1}{\delta'(y)}.
\]
\[\phi'(y)=e^{i/y}\frac{i}{y^2}\cdot \frac{\sqrt{1+y^4}}{-y^2+i}+e^{i/y}\frac{2y^3}{(-y 2+i)\sqrt{1+y^4}}+ e^{i/y} \frac{\sqrt{1+y^4}\cdot 2y}{(-y 2+i)^2}
\]
and $\delta'(y)=\displaystyle \frac{-y^2+i}{y^2e^ye^{i/y}}$.

From this, we deduce that the quotient $\frac{\phi'(y)}{\delta'(y)}$ is the sum of three terms:

The first term is
\[
\frac{e^{2i/y}i\sqrt{1+y^4}e^y}{y^2(-y^2+i)^2}
=e^ye^{2i/y}\left(\frac{iy^2\sqrt{1+\frac{1}{y^4}}}{y^4(-1+\frac{i}{y^2})^2}\right)=e^ye^{2i/y}\frac{i}{y^4}\cdot \frac{\sqrt{1+\frac{1}{y^4}}}{(-1+\frac{1}{y^2})^2}
.\]
The second term is
\[e^ye^{2i/y}\frac{2y^5}{(-y^2+i)^2\sqrt{1+y^4}}=e^ye^{2i/y}\frac{2}{y(-1+\frac{i}{y^2})^2\sqrt{\frac{1}{y^4}+1}}
.\]
The third term is
\[e^ye^{2i/y}\frac{\sqrt{1+y^4}2y^3}{(-y^2+i)^3}
.\]
 The sum of the second and the third term is
 \[e^ye^{2i/y}\left(
 \frac{2y^5}{(-y^2+i) 2\sqrt{1+y^4}}+\frac{\sqrt{1+y^4}2y^3}{(-y^2+i)^3}
 \right)\]
 \[=
 e^ye^{2i/y}\frac{2y^3}{(-y^2+i)^2}\left(\frac{y^2}{\sqrt{1+y^4}}-\frac{\sqrt{1+y^4}}{y^2-i}
 \right)
 \]
 \[= e^ye^{2i/y}\frac{2y^3}{(-y^2+i)^2}\,\frac{y^4-iy^2-1-y^4}{\sqrt{1+y^4}(y^2-i)}\]
 \[
 = e^ye^{2i/y}\frac{2(-i-\frac{1}{y^2})}{(-1+\frac{i}{y^2})^2\sqrt{\frac{1}{y^4}+1}(1-\frac{i}{y^2})}
  \]
  \[=\frac{e^ye^{2i/y}}{y^3}\cdot \left(2+\lambda(y)\right)
  \]
with $\lambda(y)\to 0$ as $y\to +\infty$.
The first term is
\[e^ye^{2i/y}\frac{i}{y^4}\frac{\sqrt{1+\frac{1}{y^4}}}{(-1+\frac{1}{y^2})^2}=\frac{e^ye^{2i/y}}{y^4}\left(i+w(y)\right)
\]
with $w(y)\to 0$ as $y\to+\infty$.

Therefore, 
\[\frac{ds^2}{d\delta^2}=e^ye^{2i/y}\left(\frac{1+w(y)}{y^4}+\frac{2+\lambda(y)}{y^3}
\right)
\]
\[=\frac{e^y}{y^3}e^{2i/y}\left(\frac{1+w(y)}{y}+2+\lambda(y)\right)
\]
\[=\frac{e^y}{y^3}e^{2i/y}\left(2+\mu(y)\right)
\]
with $\mu(y)\to 0$ as $y\to+\infty$.
Taking absolute values, we find
\[\vert \frac{ds^2}{d\delta^2}\vert =\frac{e^y}{y^3}(2+\nu(y))
\]
with $\nu(y)\to 0$ as $y\to+\infty$.

Thus, the limit of $\frac{ds^2}{d\delta^2}$ as we approach $0=\gamma(-\infty)=\delta(+\infty)$ from the curve does not exist, because its absolute value converges to $+\infty$. This implies that the curve cannot be extended beyond this point and still satisfy Definition \ref{curve-loc}. Thus, $A$ is finite.

\end{example}

\begin{example}\label{ex:5}

Let $\gamma(t)=\displaystyle  (1-x)\exp \frac{x+1}{x-1}$, $-\infty<x<1$.

Then $\gamma\subset \mathbb{R}$ and $\gamma^*=\mathbb{R}$, $A=-\infty$, $B=+\infty$, and $\lim_{x\to 1^{-}}\gamma(x)=0$. 
Thus, although the given parametrization of $\gamma$ at the one endpoint ($x\to 1^{-}$) has an essential singularity, the curve nevertheless can be continued and still satisfy Definition  \ref{curve-loc}. In fact, there is also a similar and simpler example, namely, the curve $\gamma(x)=1/x$, $-\infty <x<0$. In this case, the curve cannot be extended as $x$ converges to $0^-$ where
there is a pole , but
using arc length parametrization it can be extended as $x$ converges to $-\infty$
 where $\gamma (x)$ approaches $0$. Again, using arc length, we can extend and have the real line.
\end{example}

\begin{example}\label{ex:6}
For $\tau >0$ we consider the function
\[f(z)=(z-1)^\tau \exp\frac{z+1}{z-1}.\]
This function is holomorphic on an open set containing $\mathbb{T}-\{1\}$. Let $\gamma:(0,2\pi)\to\mathbb{C}$ be the Jordan curve $\gamma(t)=f(e^{it})$. Then, 
\[f'(z)=\left(\tau(z-1)^\tau -\frac{2}{(z-1)^{2-\tau}}\right)\cdot \exp \frac{z+1}{z-1}.\] 
Then, $0=f'(z)$, and $\tau (z-1)^\tau=\frac{2}{(z-1)^{2-\tau}}$, which implies
$\tau(z-1)^2=2$, and $z=1\pm \sqrt{\frac{2}{\tau}}$.

But $\vert 1 \pm \sqrt{\frac{2}{\tau}}\vert = 1\Rightarrow \tau=\frac{1}{2}$.
Therefore, for $\tau >\frac{1}{2}$ the function $f$ is locally conformal on an open set containing $\mathbb{T}-\{1\}= \{e^{i\theta}: 0<\theta<2\pi\}$ and $\gamma$ satisfies Definition \ref{curve-loc}. 

For $\tau >\frac{1}{2}$ we have $\lim_{t\to 0^{+}} \gamma(t)=  
\lim_{t\to 2\pi^{-}} \gamma(t)= 0$.

Furthermore, 
\[\vert \gamma'(t)\vert = \vert f'(e^{i\theta}\vert \cdot \vert ie^{it}\vert= \left\vert \tau (z-1)^\tau - \frac{2}{(z-1)^{2-\tau}}\right\vert_{z=e^{i\theta}}\]
\[\sim \vert \frac{2}{(e^{it}-1)^{2-\theta}}\vert \sim \frac{2}{t^{2-\tau}}\]
as $t\to 0^+$ ($\tau <2$) and $\vert \gamma'(t)\vert\sim \frac{2}{(2\pi-t)^{2-\tau}}$ as $t\to 2\pi^-$, since $(e^{it}-1)^\tau\to 0$ as $t\to 0^+$ or $t\to 2\pi^{-1}$ and $\vert\frac{2}{(e^{it}-1)^{2-\tau}}\vert \to +\infty$ ($\frac{1}{2}<\tau$).

Therefore, $\int_{0^{+}}\vert \gamma'(t)\vert dt=+\infty$ for $\frac{1}{2}<\tau\leq 1$.

Also, for $1<\tau$ it follows that $\int_{0^{+}}\vert \gamma'(t)\vert dt$ is finite. Similar results hold for $\int^{2\pi^{-}}\vert \gamma'(t)\vert dt$.

We conclude that for $\frac{1}{2}<\tau\leq 1$ we have $A=-\infty$ and $B=+\infty$.

For $\tau >1$ one can do a calculation similar to the one done in Example \ref{ex:4}, where we showed it is possible to have $A$ finite and $B=+\infty$.

\end{example}

\begin{example}\label{ex:7}

We take $f(z)=\left( z-( 1-\frac{1}{3})\right)^2$ and $\gamma:\mathbb{R}\to\mathbb{C}$ be defined by $\gamma(t)=f(e^{it})$. Then, $\gamma$ satisfies Definition  \ref{curve-loc}, $A=-\infty$, $B=+\infty$, $\gamma^*$ is periodic, but $\gamma$ is not a Jordan curve; it is a closed curve which intersects itself once.
Note that instead of this function, one could take any function of the form $f(z)=(z-a)^2$  where $0<a<1$, and the same property holds.

 It is known \cite{MO} that for any curve $\gamma=\phi (T)$, where $\phi$ is
holomorphic on an open set containing the unit circle $\mathbb{T}$, the complement of
$\gamma$ has finitely many components."

\end{example}

\end{document}